\documentclass[12pt]{amsart}
\usepackage{hyperref,amsmath,amssymb}

\newtheorem{theorem}{Theorem}[section]

\newtheorem{lemma}[theorem]{Lemma}

\theoremstyle{theorem}

\theoremstyle{definition}
\newtheorem{definition}[theorem]{Definition}

\theoremstyle{remark}
\newtheorem{remark}[theorem]{Remark}

\numberwithin{equation}{section}

\begin{document}
\title{On the first and second homotopy groups of manifolds with positive curvature}
\author{Richard Schoen}
\address{Department of Mathematics \\
                 University of California \\
                 Irvine, CA 92617}
\email{rschoen@math.uci.edu}
\thanks{
This work was partially supported by NSF grant DMS-2005431. }

\begin{abstract} We show that an odd dimensional closed manifold
with positive curvature cannot contain an incompressible real projective plane in the sense that there is no
map of the projective plane into the manifold which is nontrivial on both first and second homotopy groups. Another 
way to say this is that no element of the fundamental group reverses the orientation of a class in the second homotopy
group. Further we show that for a manifold of dimension divisible by four with positive curvature the fundamental group
acts trivially on the second homotopy group. The methods involve a careful study of the stability of minimal two spheres
in manifolds of positive curvature.
\end{abstract}

\maketitle

\section{Introduction} Very little is known about the topology of compact manifolds with positive sectional curvature despite the
centrality of the topic to Riemannian geometry. The Betti number bound of M. Gromov \cite{G} gives a bound on the
complexity of manifolds with nonnegative sectional curvature. The results of J. L. Synge \cite{S} using the stability of closed
geodesics imply that even dimensional manifolds of positive curvature are simply connected if and only if they are
orientable, so for such manifolds the nontrivial fundamental groups can only be $\mathbb Z_2$. Similar methods also imply that
odd dimensional manifolds with positive curvature are necessarily orientable. The restrictions for positive scalar curvature
coming from the Dirac operator also give nontrivial restrictions on manifolds of positive sectional curvature.

In this paper we add to the results of Synge which distinguish nonnegative curvature from positive curvature using the study
of the stability of minimal two spheres. In general there are stable minimal two spheres in manifolds of positive curvature, in
the even dimensional case the rational curves in complex projective spaces, and in the odd dimensional cases when 
$\pi_2\neq 0$ by the existence
theory of J. Sacks and K. Uhlenbeck \cite{SU}. We discuss the consequences of stability on the structure of the complexified normal bundle
as a holomorphic vector bundle. We show that the assumption of positive curvature induces a natural positive definite 
hermitian symmetric form on the space of holomorphic sections.

We apply this to the study of non-simply connected manifolds with positive curvature in connection with the
action of $\pi_1$ on $\pi_2$. We prove that a closed odd dimensional manifold with positive curvature cannot contain an
incompressible $RP^2$; that is, there is no map $f:RP^2\to M$ which is nontrivial on both $\pi_1$ and $\pi_2$. Another way
to describe this condition is that no deck transformation can reverse the orientation of a nontrivial class of $\pi_2$. The
additional structure which this gives us is an orientation reversing involution of $E$. Note that the minimal $S^2$ itself may
not have such a symmetry, but we are able to construct a stable minimal $S^2$ whose normal bundle has the symmetry. 
Using this involution we construct a 
complex linear form on sections which, because of the stability, has a positivity property on the space of holomorphic
sections. We then use this form in combination with the standard complex
linear form to show the instability of the minimal $S^2$. Our theorem in odd dimensions is the following.
\begin{theorem} A compact odd dimensional manifold cannot contain an incompressible $RP^2$; that is,
no element of $\pi_1$ can act by $-1$ on a nontrivial element of $\pi_2$.
\end{theorem}

Finally we obtain a result for even dimensional manifolds of dimension $4k$. It is known \cite[Section 9.2]{W} that the complex projective
space $CP^n$ has a nonorientable quotient if and only $n$ is odd, so the manifold is of dimension $2n=4k+2$. When $n$
is odd the deck transformation reverses the orientation of the standard $CP^1\subset CP^n$, and thus the quotient contains
an incompressible $RP^2$. We show that for a manifold $M^n$ of positive curvature with $n=4k$, the fundamental group
acts trivially on $\pi_2$. It follows for example that $RP^2\times S^2$ cannot have a metric with positive curvature. Note that
Synge's theorem rules out $RP^2\times RP^2$, but not $RP^2\times S^2$, More generally $RP^2\times S^{4k-2}$ does
not have a metric of positive curvature for any $k\geq 1$. Our theorem for $4k$-dimensional manifolds is the following.
\begin{theorem} Assume that $n=4k$ for a positive integer $k$, and that $M^n$ is a compact manifold with
positive curvature and with $\pi_1=\mathbb Z_2$. It follows that $\pi_1$ acts trivially on $\pi_2$.
\end{theorem}

The key to the proofs of both theorems is to use the antiholomorphic involution on $S^2$ to define a complex linear pairing on the
space of  sections of the complexified normal bundle which is associated to the hermitian pairing defined by the curvature
tensor. Using the positivity conditions on this pairing we are able to show in odd dimensions that a stable $S^2$ with normal bundle
invariant under the involution would have a holomorphic section which lies in a real subspace at each point in violation of the positivity condition.

In the $4k$-dimensional case we use Synge's theorem to show that the action of the deck transformation on the normal bundle would
necessarily preserve orientation. Using the nondegeneracy of the form we can decompose the normal bundle at a typical
point into a direct sum of two dimensional real bundles invariant under the involution and on which the involution reverses 
orientation. Since when $n=4k$ there are an odd number of such bundles and we deduce that the orientation is reversed on
the normal bundle in contradiction to Synge's theorem.

This work grew out of our recent joint work with Ailana Fraser \cite{FS} where we revisited the higher codimension stability theory
to obtain novel Bernstein-type theorems in euclidean space and quantitative results for manifolds with positive isotropic 
curvature. 

\section{Stable two spheres in manifolds with positive curvature}

In this section we consider stable minimal two spheres in manifolds with
positive curvature. Suppose we have a conformal harmonic map $f:S^2\to M$ so that $\Sigma=f(S^2)$ is a branched minimal immersion. 
We note that the tangent bundle extends across the branch points as does the normal bundle $N\Sigma$. We let $E$
denote the complexified normal bundle $E=N\Sigma\otimes\mathbb C$. On $E$ the Riemannian metric induces two 
nondegenerate pairings, the complex linear extension of the metric which we denote by $(\cdot,\cdot)$ and the other
is the positive definite Hermitian pairing denoted $\langle v, w\rangle=(v,\bar{w})$. The connection on $E$ defines a
complex structure, and $E$ becomes a holomorphic vector bundle while the complex linear pairing is holomorphic and
induces an isomorphism of $E$ with its dual bundle $E^*$. In particular $E$ has degree zero. Now by the Birkhoff-Grothendieck
theorem $E$ splits in a unique way as a direct sum of line bundles. It follows that if we write $E=P\oplus Z\oplus N$ where
$P$ is the direct sum of the positive line bundles, $Z$ the direct sum of the degree zero line bundles, and $N$ the direct
sum of the negative line bundles, then the rank of $P$ equals that of $N$ since $N=P^*$ by the self duality of $E$. 

Because we are working over $S^2$ the line bundles of nonnegative degree are spanned by holomorphic sections, and
so $P\oplus Z$ is so spanned. Now a holomorphic section $s$ of a positive line bundle has at least one zero, and so it
follows that $(s,s)$ is a holomorphic function on $S^2$ with at least one zero and therefore $(s,s)\equiv 0$. Therefore
$P$ is an isotropic bundle; that is, the complex linear pairing restricted to $P$ vanishes identically.  Note that a holomorphic
section $s$ of a degree zero bundle has no zeroes, so the function $(s,s)$ is constant, but not necessarily zero. 

We choose $x,y$ to be positively oriented isothermal coordinates on $S^2$.
The stability condition for area implies that for any smooth section $s$ of $E$ we have
\begin{equation}\label{stability_0} \int_{S^2}R^M(f_x+if_y,s,f_x-if_y,\bar{s})\ dx dy\leq 4\int_{S^2}\|\nabla_{\bar{z}}^\perp s\|^2\ dx dy
\end{equation}
where $x,y$ are local oriented isothermal coordinates, and $f_x,f_y$ denote the associated orthogonal tangent basis of $\Sigma=f(S^2)$. 
It follows that if $s$ is holomorphic we have 
\[ \int_{S^2}R^M(f_x+if_y,s,f_x-if_y,\bar{s})\ dx dy\leq 0.
\]
In order to interpret this condition, we use the notation $v_1=f_x$, $v_2=f_y$, and $s=v_3+iv_4$. We then
note that $(v_1,v_2,v_3,v_4)$ is a $4$-tuple of vectors which span a space of at least three dimensions at all but a finite set of points. This is because of the isothermal condition and the fact that
$s$ is a nonzero section of $E$. We can then expand the complex sectional term out to obtain
\[ R^M(f_x+if_y,s,f_x-if_y,\bar{s})=R_{1313}+R_{2323}+R_{1414}+R_{2424}-2R_{1234}.
\]
Note that each of the first four terms is non-negative and the sum is strictly positive at all but a finite set of points since $s=v_3+iv_4$ is a nonzero holomorphic section of $E$. We may therefore conclude that for a nonzero holomorphic section of $E$ we have
\[ \int_{S^2} R_{1234}\ dxdy>0.
\]
If we choose a positively oriented local orthonormal tangent basis $e_1,e_2$ we may rewrite this condition
\[ \int_{S^2}\ R(e_1,e_2,v_3,v_4)\ da>0,
\]
where $da$ denotes the (unoriented) area measure of $\Sigma$. Writing this back in terms of $s=v_3+iv_4$ we have
\begin{equation}\label{pos_form} \sqrt{-1}\int_{S^2} R(e_1,e_2,s,\bar{s})\ da>0.
\end{equation}

While stability is an integral condition, it motivates the definition of a Hermitian form on the bundle $E$ given for any point $x\in S^2$ by
\[ \langle v,w\rangle_x=\sqrt{-1}R(e_1,e_2,v,\bar{w})
\]
where $v,w\in E_x$ and $e_1,e_2$ is any chosen oriented orthonormal basis for the tangent space at $x$.
Note that $\langle\cdot,\cdot\rangle_x$ is linear in the first slot and conjugate linear in the second, and from the symmetries of the curvature tensor 
we have $\langle v,w\rangle_x=\overline{\langle w,v\rangle_x}$. Thus $\langle\cdot,\cdot\rangle_x$
defines a Hermitian symmetric form at each point on the fiber $E_x$. Notice also that if $v$ is a complex number times
a real vector we have $\langle v,v\rangle_x=0$. We may then rewrite \ref{stability_0} for a section $s=v_1+iv_2$ as 
\begin{equation}\label{stability} 
\int_{S^2} (R_{1313}+R_{1414}+R_{2323}+R_{2424}-\langle s,s\rangle_x)\ da\leq 4\int_{S^2}\|\nabla^\perp_{\bar{z}}s\|^2\ dxdy.
\end{equation}

We can then define the integrated form $\langle\cdot,\cdot\rangle_0$ on the space of sections of $E$ by
\[ \langle s_1,s_2\rangle_0=\int_{S^2}\langle s_1,s_2\rangle_x\ da.
\]
We also point out the following antisymmetry property which follows directly from the definition 
\[
\langle s_1,\bar{s}_2\rangle=-\langle s_2,\bar{s}_1\rangle.
\]
Now if we let $H(E)$ denote the space of holomorphic sections of $E$ we have the following lemma.

\begin{lemma}\label{omega-plus} If $M^n$ is a manifold with positive curvature, $f:S^2\to M$ a non-constant harmonic mapping
which is stable, then the form $\langle\cdot,\cdot\rangle_0$ is positive definite on $H(E)$. 
\end{lemma}

\section{Representing incompressible projective planes}

Let $M^n$ be a compact Riemannian manifold with $n\geq 3$ and with finite fundamental group. We will say that $M$
contains an {\it incompressible $RP^2$} if there is a smooth map $f:RP^2\to M$ which is nontrivial on both $\pi_1$
and $\pi_2$. Note that if we lift $f$ to universal coverings we produce a map $\tilde{f}:S^2\to\tilde{M}$ and a deck transformation
$\gamma$ so that $\tilde{f}\circ\gamma_0=\gamma\circ \tilde{f}$ where $\gamma_0$ is the deck transformation for $RP^2$,
the antipodal map of the unit $S^2$. The deck transformation $\gamma$ then maps the class $[\tilde{f}]\in \pi_2(M)$ to its 
negative. 

Our goal is to represent an incompressible $RP^2$ by a least area (or energy) map, and then to study the structure of such stable
minimal immersions. If we minimize the energy among maps $\tilde{f}$ with the condition that $\tilde{f}\circ\gamma_0=\gamma\circ \tilde{f}$,
we will get a minimal $S^2$, but it will not necessarily be stable for all variations, but only for those which preserve the symmetry. For
example the standard embedding of $RP^2$ into $RP^n$ for $n\geq 3$ is injective on $\pi_1$ but not on $\pi_2$, and its lift to $S^n$
is stable for variations which preserve the symmetry, but is unstable for general variations. 

We have the following existence theorem for stable two spheres which uses the refinement of the work of Sacks-Uhlenbeck \cite{SU}
done by Meeks-Yau \cite[Theorem 5]{MY}.

\begin{theorem}\label{sphere-existence}  Let $f$ be an incompressible $RP^2$ and let $\gamma$ be the associated deck transformation
which reverses the $\pi_2$ class defined by a lift of $f$. There are a pair of not necessarily distinct stable $S^2$'s parametrized 
conformally by maps $f_1$ and $f_2$ from $S^2$ into $\tilde{M}$ so that $\gamma\circ f_1=f_2\circ\gamma_0$ and so that
$\gamma\circ f_1\circ\gamma_0$ is homotopic to $f_1$. 
\end{theorem}

\begin{proof} If we assume that we can find a least energy map $f_1:S^2\to \tilde{M}$ in the homotopy class of $\tilde{f}$,
then we can define $f_2\equiv \gamma\circ f_1\circ\gamma_0$, and we have $f_2\circ\gamma_0=\gamma\circ f_1$ as
claimed. Note that the maps $f_1$ and $f_2$ are homotopic, but may or may not coincide.

While it is not necessarily true that we can find an energy minimizing map in the homotopy class of $\tilde{f}$, the work of
Sacks-Uhlenbeck \cite{SU} tells us that we can find stable harmonic maps $v_1,\ldots, v_k$ from $S^2$ into $\tilde{M}$ which
are nontrivial in $\pi_2$ and such that the the homotopy classes satisfy $[\tilde{f}]=[v_1]+\ldots +[v_k]$. The work of Meeks-Yau \cite{MY}
tells us that we can take the harmonic maps $v_j$ so that they generate $\pi_2$ as a module over $\pi_1$. Since we have
$[\gamma\circ\tilde{f}]=-[\tilde{f}]$ it follows that we may assume that for each $j$ we have $[\gamma\circ v_j]=-[v_j]$. It is clear
that $[\gamma\circ v_j\circ\gamma_0]=-[\gamma\circ v_j]$, and thus $[\gamma\circ v_j\circ\gamma_0]=[v_j]$. We can now
define $f_1=v_1$ and $f_2=\gamma\circ f_1\circ\gamma_0$, and both are stable minimal spheres that satisfy the required
relation. 
\end{proof}

If the maps $f_1$ and $f_2$ are the same, we would have a stable minimal sphere in $\tilde{M}$ defined by $f_1$ with $f_1\circ\gamma_0=\gamma\circ f_1$, and we could define a map $f:RP^2\to M$ since the projection $\pi:\tilde{M}\to M$ 
identifies the points $x$
and $\gamma(x)$. Thus if we let $E$ be the complexified normal bundle of the branched immersion $f$, we could lift $E$ to
$S^2$ to get a bundle with the symmetry $\gamma_0^*E=E$.

In case $f_1\neq f_2$ we cannot define a map from $RP^2$ to $M$ which is covered by $f_1$, but we can still define a
bundle $E$ with the same symmetry. The idea is that we have two minimal spheres which are interchanged by $\gamma$.
If we parametrize one of them by a map $f:S^2\to M$, then for any $x\in S^2$ we can consider the pair of points 
$f(x)$ and $f\circ\gamma_0(x)$. These are mapped by $\gamma$ to a pair of points on the other sphere. We will show that the
parametrization can be chosen so that the pair of image points is $\gamma\circ f\circ\gamma_0(x)$ and $\gamma\circ f(x)$; in
other words we can parametrize by a map $f$ with $\gamma\circ f\circ\gamma_0=\gamma\circ f$.
Since these points are identified in the quotient, it follows that the pulled back normal bundle is invariant under $\gamma_0$.

To see this rigorously we consider the maps $\pi\circ f_1$ and $\pi\circ f_2$. where $\pi:\tilde{M}\to M$ is
the covering projection These maps have
the same image $\Sigma$ which is a stable $S^2$ in $M$ and has complexified normal bundle $N\Sigma\otimes \mathbb C$. The maps
are homotopic since the maps $f_1$ and $f_2=\gamma\circ f_1\circ\gamma_0$ are homotopic. Therefore there is an oriented conformal diffeomorphism $\phi$ 
of $S^2$ so that $\pi\circ f_1=\pi\circ f_2\circ\phi$. We clearly also have $\pi\circ f_2=\pi\circ f_1\circ\gamma_0$ and so 
$\pi\circ f_2=\pi\circ f_2\circ\phi\circ\gamma_0$. Now
$\phi$ is in the identity component of the conformal diffeomorphism group of $S^2$, so $\phi\circ\gamma_0$ is in the component of
the space of orientation reversing conformal diffeomorphisms containing $\gamma_0$. It follows that there is an orientation
preserving conformal diffeomorphism $\psi$ of $S^2$ so that $\phi\circ\gamma_0=\psi\circ\gamma_0\circ\psi^{-1}$. Thus
we have $\pi\circ f_2=\pi\circ f_2\circ\psi\circ\gamma_0\circ\psi^{-1}$, and so if we
define $f:S^2\to\tilde{M}$ by $f=f_2\circ\psi$, we have $\pi\circ f=\pi\circ f\circ\gamma_0$ and $\pi\circ f$ is a stable branched
immersion. We then define 
$E=(\pi\circ f)^*N\Sigma\otimes\mathbb C$ and we have $\gamma_0^*E=E$ as claimed. In the next section we will use this 
symmetry in an important way to obtain our main results.

\section{The action of $\pi_1$ on $\pi_2$ for manifolds of positive curvature}

We now assume that $M$ is a compact manifold of dimension $n\geq 4$ with a metric of positive sectional curvature. Recall
from Synge's theorems that if $n$ is even then $\pi_1$ is trivial if $M$ is orientable, while $\pi_1=\mathbb Z_2$ if $M$ is
not orientable. Also, if $n$ is odd, then $M$ is necessarily orientable. 

We first digress to discuss some linear algebra preliminaries. We assume that $\mathbb C^n$ is equipped with a symmetric complex linear
pairing $(\cdot,\cdot)_0$. Assume the pairing is pure imaginary in that $(\bar{v},\bar{w})_0=-\overline{(v,w)_0}$. Suppose also that there is a linear map $T:\mathbb C^n\to \mathbb C^n$ such that $T^2=I$, $\overline{Tv}=T\bar{v}$, and that 
$(Tv,Tw)_0=-(v,w)_0$ for all $v,w\in\mathbb C^n$. 

\begin{definition} We will say that the form $(\cdot,\cdot)_0$ is {\it T-positive} on a subspace $V^k\subset\mathbb C^n$ if $T(V)=\bar{V}$
and for all $v\in V$, we have $(v,T\bar{v})_0\geq 0$ and strictly positive unless $v=0$.
\end{definition}
Note that $T(V)=\bar{V}$ if and only if $T(\bar{V})=V$ since $T^2=I$.

We remark that the form $\langle v,w\rangle_0\equiv (v,T\bar{w})_0$ is a sesquilinear form in the sense that it is linear in the first slot,
conjugate linear in the second and $\langle v,w\rangle_0=\overline{\langle w,v\rangle_0}$. To see the last property, note that from the
properties of $T$ we have $(v,T\bar{w})_0=-(\bar{w},Tv)_0=\overline{(w,T\bar{v})_0}$, and this implies 
$\langle v,w\rangle_0=\overline{\langle w,v\rangle_0}$.

Note that if the form is T-positive on $V$ then the complex linear pairing $(\cdot,\cdot)_0$ is clearly nondegenerate on $V$. We will need the following elementary result.

\begin{lemma}\label{lin_alg_1} With the notation as above, let $V^k$ be a subspace on which the form is T-positive and with $T(V)=\bar{V}$. 
The form is then nondegenerate on $V\oplus\bar{V}$ and there is a basis of the form $v_1,Tv_1,v_2,Tv_2,\ldots, v_k,Tv_k$ with 
$Tv_j=\bar{v}_j$ and which
is orthonormal in the sense that $(v_j,v_k)_0=\delta_{jk}$, $(v_j,Tv_k)_0=0$, and $(Tv_j,Tv_k)_0=-\delta_{jk}$. 
\end{lemma}

\begin{proof} Let $w_1$ be a unit vector in $V$. If $w_1=-T\bar{w}_1$ replace $w_1$ by $iw_1$, so we may assume that 
$w_1\neq -T\bar{w}_1$. Now let 
\[ v_1=\frac{w_1+T\bar{w}_1}{\sqrt{(w_1+T\bar{w}_1,w_1+T\bar{w}_1)_0}}
\]
where we note the term under the square root is positive by T-positivity. We thus have $v_1$ is a unit element of $V$ with $Tv_1=\bar{v}_1$.
We note that $(v_1,v_1)_0=1$, $(Tv_1,Tv_1)_0=-1$, and $(v_1,Tv_1)_0=0$ (the third is because the form $(v,Tw)_0$ is anti-symmetric).
We now consider the orthogonal complement of the subspace spanned by $v_1,Tv_1$ and we observe that it is direct sum of
$V_1=v_1^\perp\subset V$ with its conjugate, and $T(V_1)=\bar{V}_1$. We then repeat the argument to produce $v_2\in V_1$ with
$Tv_2=\bar{v}_2$. Since $V_1\oplus \bar{V}_1$ is orthogonal to the space spanned by $v_1,Tv_1$ it follows that $v_2$ and $Tv_2$ 
are both orthogonal to both $v_1$ and $Tv_1$. We may then repeat this argument to produce the basis claimed. In particular this
argument shows that the form is nondegenerate on $V\oplus\bar{V}$ competing the proof of Lemma \ref{lin_alg_1}.
\end{proof}

We prove another linear algebra result which will be used in the odd dimensional case.
\begin{lemma}\label{lin_alg_2} With the notation as above and for a positive integer $l$, if the form is T-positive on $V^{l+1}$ and 
$V_1\subset V$ is a subspace of dimension $l$ with $T(V_1)=\bar{V}_1$ and $T(V)=\bar{V}$, then the form is nondegenerate on $V\oplus \bar{V}_1$ and there is an element $v_0\in V$ orthogonal to $V_1\oplus\bar{V}_1$. 
\end{lemma}

\begin{proof} We may use Lemma \ref{lin_alg_1} to show that the form is nondegenerate on $V_1\oplus \bar{V}_1$ and may
be diagonalized in an orthonormal basis $v_1,Tv_1,\ldots, v_l,Tv_l$. Now we can repeat this to extend the basis to one of
$V\oplus\bar{V}$ by adding $v_{l+1},Tv_{l+1}$. The element $v_{l+1}\in V$ thus satisfies the requirements for $v_0$ claimed
in the result. This completes the proof of Lemma \ref{lin_alg_2}.
\end{proof}

We first set up the notation for our main results. Assume we have the complexified normal bundle $E$ from a branched
minimal immersion in a Riemannian manifold $M$ with its metric and connection. We assume also that the bundle is
invariant under the antipodal map $\gamma_0:S^2\to S^2$ so that $\gamma_0^*E=E$ and $\gamma_0$ preserves the metric
and connection. The vector spaces we will consider will be subspaces of the space of smooth sections $\Gamma(E)$. We define 
a map $T:\Gamma(E)\to\Gamma(E)$ by setting $Ts=s\circ\gamma_0$; that is, $Ts$ is the pullback of $s$ under $\gamma_0$.
Note that $T^2=I$ and $\overline{Ts}=T\bar{s}$.
We then define the pairing $(\cdot,\cdot)_0$ on $\Gamma(E)$ by
\[ (s_1,s_2)_0=\sqrt{-1}\int_{S^2}R(e_1,e_2,s_1, Ts_2)\ da.
\]
The pairing is symmetric as we can see by doing a change of variables using $\gamma_0$. This operation changes $s_1$ and
$s_2$ to $s_1\circ\gamma_0$ and $s_2\circ\gamma_0$, so we get a sign reversal upon interchange of the second two slots,
but $\gamma_0$ also reverses the orientation of the tangent plane, so $e_1$ and $e_2$ also get reversed and we end up showing that
$(s_1,s_2)_0=(s_2,s_1)_0$. The pairing is clearly complex linear and satisfies $(Ts_1,Ts_2)_0=-(s_2,s_1)_0=-(s_1,s_2)_0$ and 
$(\bar{s}_1,\bar{s}_2)_0=-\overline{(s_1,s_2)_0}$. More precisely it is symmetric because
\[ \int_{S^2}R(e_1,e_2,s_1,s_2\circ\gamma_0)\ da=\int_{S^2}R(e_1\circ\gamma_0,e_2\circ\gamma_0,s_1\circ\gamma_0,s_2)\ da
\]
since $\gamma_0$ is an isometry which preserves the curvature tensor. Since $\gamma_0$ reverses orientation, the orthonormal 
basis $e_1\circ\gamma_0,e_2\circ\gamma_0$ is negatively oriented, but the third and fourth arguments may also be reversed so
that we have $(s_1,s_2)_0=(s_2,s_1)_0$.

Notice that the hermitian pairing defined in the first section of the paper is given by
\[ \langle s_1,s_2\rangle_0=(s_1,T\bar{s}_2)_0,
\]
so it follows from (\ref{pos_form}) that if $M$ has positive curvature, the paring is T-positive on H(E), the space of holomorphic sections of $E$.

Our first main result of the section is that for $n$ odd there can be no incompressible $RP^2$ in a manifold $M^n$ of positive curvature.
\begin{theorem}\label{odd_proj_plane} A compact odd dimensional manifold cannot contain an incompressible $RP^2$; that is,
no element of $\pi_1$ can act by $-1$ on a nontrivial element of $\pi_2$.
\end{theorem}

\begin{proof} In the previous section we showed that if $M$ has an incompressible $RP^2$ then there is a stable minimal
immersion $f:S^2\to M$ so that the complexified normal bundle $E$ is invariant under $\gamma_0$. Using the normal connection
we can make $E$ into a holomorphic vector bundle, and using the Birkhoff-Grothendieck theorem we can split $E$ as a direct
sum of holomorphic line bundles $E=L_1\oplus \ldots \oplus L_{n-2}$. If we have a local holomorphic section $s$, then $s\circ\gamma_0$
is a local anti-holomorphic section, and thus $\bar{s}\circ\gamma_0$ is a local holomorphic section. It follows that if $L_j$ is
a holomorphic line sub-bundle of $E$ then so is $\gamma_0^*\bar{L_j}$. We also note that the degree of $\gamma_0^*\bar{L_j}$ is the same
as that of $L_j$ because if $s$ is a meromorphic section of $L_j$, then the section $\bar{s}\circ\gamma_0$ is a meromorphic section of 
$\gamma_0^*\bar{L_j}$ with the same number of zeroes and poles. By the uniqueness of the splitting it follows that if we write $E=P\oplus Z\oplus N$ where $P$, $Z$, $N$ denote the direct sum of the positive, zero, negative bundles, then $\gamma_0^*P=\bar{P}$, 
$\gamma_0^*Z=\bar{Z}$, and $\gamma_0^*N=\bar{N}$. Since $E$ is self dual with isomorphism given by the nondegenerate complex linear pairing defined by the Riemannian metric, we know that $N=P^*$ and so they are of the same rank. Since $n$ is odd $Z$ must have odd rank. Since
we are working over $S^2$, the bundle $Z$ is a holomorphically trivial bundle. 

Let the rank of $Z$ be $2q+1$ for a nonnegative integer $q$. If $q>0$, we want to choose a $q+1$ dimensional subspace of the
holomorphic sections $H(Z)$ in such a way that it decomposes as a direct sum of a subspace of isotropic sections (for the usual
complex linear pairing) together with a one dimensional subspace which is nondegenerate for the same pairing. Note that
if $s_1,s_2\in H(Z)$ then $(s_1,s_2)$ is a holomorphic function on $S^2$ which is thus a constant. Thus the pairing $(s_1,s_2)$
is a complex valued pairing. We show that this pairing is nondegenerate on $H(Z)$. To see this note that the same argument
shows that for any $x\in S^2$ the space $P_x$ is an isotropic subspace. This is because each positive line bundle has a holomorphic
section with at least one zero, and therefore the pairing between any two of these sections is zero. Similarly the space $Z_x$
is orthogonal to $P_x$ with respect to $(\cdot,\cdot)$. 

To see that $(\cdot,\cdot)$ is nondegenerate on $H(Z)$, we suppose it had a null space $K\subset H(Z)$ of dimension $r$. We could then
choose a complementing subspace $W$ so that $H(Z)=K\oplus W$. Now the pairing must be nondegenerate on $W$, and so we
can choose an orthonormal basis $s_1,\ldots, s_t$ where $t=2q+1-r$ for $W$. We can then consider the sections of the form
$s_1+is_2,s_3+is_4,\ldots$ to produce a space of isotropic sections with dimension $[t/2]$, the integer part of $t/2$. If we combine
these with the space $K$ of sections, we get a space of isotropic sections of $H(Z)$ of dimension $r+[(2q+1-r)/2]$ and this
number is larger than $[(2q+1)/2]$ if $r>0$. If we combine this isotropic subspace at a point $x$ with the subspace $P_x$ we 
get an isotropic subspace of $E$ with dimension greater than $[n/2]$. Since the pairing is nondegenerate on $E$, any
isotropic subspace $V$ of $E_x$ satisfies $V\cap\bar{V}=\{0\}$, and therefore the dimension of $V$ can be at most $[n/2]$.
Therefore we have $r=0$ and the pairing is nondegenerate on $H(Z)$.

We also observe that $T(H(Z))=\overline{H(Z)}$, so we can choose an orthonormal basis $s_1,\ldots s_{2q+1}$ for $H(Z)$ with 
respect to $(\cdot,\cdot)$ so that $Ts_j=\bar{s}_j$ and from it construct a $q$ dimensional
isotropic subspace $W_1\subset H(Z)$ (recall $\dim H(Z)=2q+1$) with $T(W_1)=\overline{W}_1$. The remaining section in the orthonormal basis then spans a one dimensional nondegenerate subspace $W_2$ of $H(Z)$ with $T(W_2)=\overline{W}_2$. 

We now consider the positive bundle $P$ of rank $p$, and we note that $T(H(P))=\overline{H(P)}$. More precisely we have positive line bundles $L_j$ so that $P=L_1\oplus\ldots\oplus L_p$ and there is a permutation $\sigma$ of order $2$ so that 
$T(L_j)=\overline{L}_{\sigma(j)}$. If we choose nonzero holomorphic sections $s_j$ of $L_j$, then $s_1,\ldots,s_p$ form a basis for $P_x$ at all but a finite set of points of $S^2$; that is, away from the points where one of the $s_j$ is zero. We construct our sections to respect
the action of $T$ as follows: if $\sigma(j)=j$, we let $s$ be a nonzero section and by replacing $s$ with $is$ if necessary we
may assume that $T\bar{s}\neq -s$. We then let $s_j=s+T\bar{s}$ and we have $Ts_j=\bar{s}_j$ (note that in this case the degree of $L_j$
is even since the zeroes of $s_j$ come in antipodal pairs). If $\sigma(j)\neq j$, we choose any nonzero holomorphic section
of $L_j$ and take $s_{\sigma(j)}=T\bar{s}_j$. By this choice we have $T(L_j\oplus L_{\sigma(j)})=\overline{L_j\oplus L_{\sigma(j)}}$. We let $W_3$ be the subspace of $H(E)$ spanned by $s_1,\ldots,s_p$, and by construction we have $T(W_3)=\overline{W}_3$ since 
any member of $W_3$ can be written $s=\sum_{j=1}^p a_js_j$ for complex numbers $a_j$, so its image under $T$ is given by
$\sum_{j=1}^pa_{\sigma(j)} \bar{s}_j$, and this is the conjugate of $\sum_{j=1}^p\bar{a}_{\sigma(j)}s_j\in W_3$.

Finally we form a space $V$ of dimension $k+1$ where $n=2k+1$ by letting 
\[ V=W_1\oplus W_2\oplus W_3. 
\]
From the construction we see that, with respect to the form $(\cdot,\cdot)$, $W_1\oplus W_3$ is isotropic and $W_2$ is nondegenerate. 
We now let $V_1$ be the codimension one subspace of $V$ given by
\[ V_1=W_1\oplus W_3.
\]
We note that away from a finite set of points of $S^2$, the space $(V\oplus \bar{V}_1)_x$ spanned by the members of $
V\oplus \bar{V}_1$ evaluated at $x$ is $E_x$.

We now see that the pair $V_1\subset V$ satisfy the hypotheses of Lemma \ref{lin_alg_2} with respect to $T$ and $(\cdot,\cdot)_0$,
so there exists a nonzero $s\in V$ which is orthogonal with respect to $(\cdot,\cdot)_0$ to $V_1\oplus \bar{V}_1$. Thus if we
go to a point $x\in S^2$ where $s$ is nonzero and where $E_x=(V\oplus \bar{V}_1)_x$, we see that the value $s(x)$ is
orthogonal to a real subspace of codimension one with respect to a nondegenerate $(v,w)_x$ defined by $(s_v,s_w)_0$ where $s_v,s_w$
are the unique members of $V\oplus \bar{V}_1$ with $s_v(x)=v,\ s_w(x)=w$. The conditions satisfies by $(\cdot,\cdot)_0$ then
imply that this subspace is a real one dimensional subspace and thus $s(x)$ is a complex number times a real vector. This
property then holds at all points of $S^2$. This is a contradiction since from \ref{pos_form} we have
\[ \sqrt{-1}\int_{S^2}R(e_1,e_2,s,\bar{s})\ da>0,
\]
but the integrand is identically zero if $s$ is a real section.  

To make this argument precise, we let $(\cdot,\cdot)_x$ denote the pairing induced from $(\cdot,\cdot)_0$ by the 
evaluation map from $V\oplus\bar{V}_1$ to $E_x$. We let $V_0$ be the linear span of the element $s\in V$ which
is orthogonal to $V_1\oplus\bar{V}_1$ with respect to $(\cdot,\cdot)_0$ and we note if $\bar{s}\in V_0$ then $s$ would
be real which as shown above is not possible. Thus $\bar{V}_0\cap V_0=\{0\}$ and from Lemma \ref{lin_alg_2}
we have $V_0\oplus\bar{V}_0$ orthogonal to $V_1\oplus\bar{V}_1$. It follows that the evaluation map is an isometry
from the orthogonal complement of $V_0\oplus \bar{V}_0$ to the orthogonal complement of $s(x)\in E_x$. It follows that
the image under the evaluation map of $V_0\oplus\bar{V}_0$ is the linear span of $s(x)$ which is a one dimensional
subspace of $E_x$. This subspace is real since it is the orthogonal complement of a real subspace. Thus we get a contradiction which 
shows that $M$ cannot contain an incompressible $RP^2$. This completes the proof of Theorem \ref{odd_proj_plane}.
\end{proof}

\begin{remark} An example of a manifold for which our theorem applies which is not covered by Synge's theorem is an
orientable $\mathbb Z_2$ quotient of $S^2\times S^3$ given by an isometry of the form $\gamma(x,y)=(-x,Ay)$ where
$A$ is an orientation reversing isometry of order $2$ of $S^3$. Note that $\gamma$ acts freely even though $A$ has
fixed points.
\end{remark}

We now come to the case of even dimensional manifolds. In this case Synge's theorem implies that a nontrivial fundamental group
can only be $\mathbb Z_2$ and this if and only if $M$ is not orientable. In the case of $CP^n$, a non-orientable quotient
exists only in the case of $n$ odd. It is interesting to ask if manifolds of positive curvature share a property of this type. We 
prove the following result.

\begin{theorem}\label{even_proj_plane} Assume that $n=4k$ for a positive integer $k$, and that $M^n$ is a compact manifold with
positive curvature and with $\pi_1=\mathbb Z_2$. It follows that $\pi_1$ acts trivially on $\pi_2$.
\end{theorem}

\begin{proof} Let $\gamma$ be the deck transformation acting on $\tilde{M}$ so that $\gamma^2=I$ and $M=\tilde{M}/<\gamma>$.
If there is an element $\alpha\in \pi_2$ such that $\gamma^*\alpha\neq \alpha$, then $\beta=\alpha-\gamma^*\alpha$ satisfies
$\gamma^*\beta=-\beta$, and we have an incompressible $RP^2$ in $M$. 

By the existence theory we then have a stable minimal sphere parametrized by $f:S^2\to M$ with the property that the complexified
normal bundle $E$ is invariant under $\gamma_0$. Because $\gamma_0$ reverses orientation on $f^*TM=f^*T\Sigma\oplus N\Sigma$
and $\gamma_0$ reverses the orientation of $T\Sigma$ it follows that $\gamma_0$ preserves the orientation of the normal bundle
$N\Sigma$. 

Now we let $E=N\Sigma\otimes \mathbb C$, and split $E$ as a direct sum of holomorphic line bundles. Combining the bundles of
positive, zero, and negative degree we have $E=P\oplus Z\oplus N$. If $Z$ was of rank $0$, then $P$ would be a $2k-1$ dimensional
subbundle of $E$. Since we are working over $S^2$, $P$
is spanned by holomorphic sections each of which has at least one zero. Since $(s,s)$ is a constant for any holomorphic section $s$
we see that $P_x$ is an isotropic subspace of $E_x$ for all $x\in S^2$. We may then use $P$ to define an orthogonal complex
structure on $E$. We do this by defining $J=iI$ on $P$ and $J=-iI$ on $\bar{P}$. This in turn defines an orientation on $N\Sigma$
by the requirement that a basis of the form $e_1,Je_1,\ldots, e_{2k-1},Je_{2k-1}$ be positively oriented. 

We showed in the proof of Theorem \ref{odd_proj_plane} that $\gamma_0^*P=\bar{P}$, and, in fact, if we write 
$P=L_1\oplus\ldots\oplus L_{2k-1}$, there is a permutation $\sigma$ of $\{1,2,\ldots, 2k-1\}$ of order two such that 
$\gamma_0^*L_j=\bar{L}_{\sigma(j)}$. It follows that $\gamma_0$ acts by complex conjugation on the bundle 
$L_j\oplus \bar{L}_{\sigma(j)}$. Now we may write 
\[ E=P\oplus\bar{P}=[L_1\oplus\bar{L}_{\sigma(1)}]\oplus\ldots\oplus [L_{2k-1}\oplus\bar{L}_{\sigma(2k-1)}].
\]
We show that $\gamma_0$ reverses the orientation of $N\Sigma$ by noting that if we choose a local oriented orthonormal basis
$v_1,w_1,v_2,w_2\ldots, v_{2k-1},w_{2k-1}$ so that $v_j+iw_j$ spans $L_j$, then we have
\[ v_1\wedge w_1\wedge\ldots\wedge v_{2k-1}\wedge w_{2k-1}=(i/2)^{2k-1}(v_1+iw_1)\wedge
(v_1-iw_1)\wedge\ldots\wedge (v_{2k-1}+iw_{2k-1}\wedge (v_{2k-1}-iw_{2k-1}).
\]
We see that the action of $\gamma_0$ switches $v_j+iw_j$ with $v_{\sigma(j)}-iw_{\sigma(j)}$ Since there are an odd number
of $j$, it follows that $\gamma_0$ reverses the orientation of $E$, a contradiction which completes the proof of
 Theorem \ref{even_proj_plane} in case $Z$ is trivial.
 
 If $Z$ is nontrivial, then the rank $p$ of $P$ is less than $2k-1$, and the rank of $Z$ is $2(2k-1-p)$. We can proceed as in
 the proof of Theorem \ref{odd_proj_plane} and since the complex linear pairing $(\cdot,\cdot)$ is nondegenerate on $Z$
 we can construct a subspace $Z_0$ consisting of isotropic holomorphic of dimension $2k-1-p$. We now consider $P_0\equiv P\oplus Z_0$
 which is an isotropic holomorphic subbundle of $E$ of rank $2k-1$. We see that $E=P_0\oplus\bar{P}_0$ and that by construction
 $P_0$ is a direct sum of nonnegative line bundles. Therefore we can find a space $V$ with dimension $2k-1$ of holomorphic sections
 with $T(V)=\overline{V}$ and so that the evaluation map is a linear isomorphism from $V\oplus\bar{V}$ to $E_x$ at all but a finite set of points of $S^2$. We may now apply Lemma \ref{lin_alg_1}
 to construct an orthonormal basis (with respect to $(\cdot,\cdot)_0$) of $V\oplus\bar{V}$ of the form $s_1,Ts_1,\ldots, s_{2k-1},Ts_{2k-1}$ where $Ts_j=\bar{s}_j$ for $j=1.\ldots,2k-1$. Now we let $x\in S^2$ be a point at which the evaluation map from $V\oplus\bar{V}$ is an isomorphism
 to $E_x$, and we let $(\cdot,\cdot)_x$ denote the pairing on $E_x$ induced by this isomorphism. We then let $v_j=s_j(x)$,
 and we see that letting $V_j$ be the linear span of $v_j$, we can write 
 $E_x=(V_1\oplus TV_1)\oplus\ldots\oplus (V_{2k-1}\oplus TV_{2k-1})$ and each of the two dimensional subspaces $V_j\oplus TV_j$ 
 is a real subspace. If we let $v_j=u+iw$ for real vectors $u,w$ then because $v_j$ is isotropic and nonzero, the vectors $u$ and
 $w$ are orthogonal and of equal nonzero length with respect to $(\cdot,\cdot)$. Now the action of $\gamma_0$ is given by $T$
 and we see that $Tv_j=\bar{v_j}$ means that $Tu=u$ and $Tw=-w$. Therefore the action of $\gamma_0$ reverses orientation
 on each of the real summands $V_j\oplus TV_j$, and since there are an odd number of them, it reverses orientation on $E_x$. 
 This contradiction completes the proof of Theorem \ref{even_proj_plane}
\end{proof}

\begin{remark} The complex projective space $CP^n$ has dimension $2n$, and it has a nonorientable quotient for $n$ odd
with the orientation reversing deck transformation also reversing the orientation of the standard $CP^1\subset CP^n$ 
(see \cite[Section 9.2]{W}).
Thus this quotient has an incompressible $RP^2$. 
Theorem \ref{even_proj_plane} shows that manifolds of positive curvature do not have nonorientable quotients in dimension 
$4k$ by a deck transformation acting nontrivially on $\pi_2$, sharing this property with $CP^n$. Note that for $CP^n$ the Euler 
characteristic is $n+1$, and this is odd when $n$ is even.
This also implies that there can be no nonorientable quotient for $n$ even. We know little about the Euler characteristic of
a manifold of positive curvature, so this argument does not extend.
\end{remark}

\begin{remark} Theorem \ref{even_proj_plane} implies that $RP^2\times S^{4k-2}$ cannot have a metric of positive curvature
for any positive integer $k$. This conclusion does not seem to follow from geodesic arguments.
\end{remark}

\bibliographystyle{plain}

\begin{thebibliography}{BFNT}

\bibitem{G} M. Gromov, Curvature, diameter and Betti numbers,
                 {\em Comment. Math. Helv.} {\bf 56} (1981), no. 1 179--195.
\bibitem{FS} A.Fraser and R. Schoen, Stability and largeness properties of minimal surfaces in higher codimension, 
               arXiv:2303.07423.
\bibitem{MY}  W. Meeks and S.~T. Yau, Topology of three-dimensional manifolds and the embedding problems in
              minimal surface theory, 
              {\em Ann. of Math. (2)} {\bf 112} (1980), no. 3, 441--484.    
\bibitem{SU} J. Sacks and K. Uhlenbeck, , The existence of minimal immersions of $2$-spheres, 
              {\em Ann. of Math. (2)} {\bf 113} no. 1 (1981), 1--24.                                       
\bibitem{S} J.~L. Synge, On the connectivity of spaces of positive curvature, 
             {\em Quarterly Journ. of Math.}, Oxford series {\bf 7} (1936), no. 1, 316--320.
\bibitem{W} J.~A. Wolf, {\em Spaces of Constant Curvature}, McGraw-Hill, New York-Sydney (1967).
                        
\end{thebibliography}

\end{document}